\documentstyle[a4,amstex,12pt]{article}
     \newcommand{\lrw}{\longrightarrow}

\newcommand{\coz}{{\operatorname{coz}}}
\newcommand{\supp}{{\operatorname{supp}}}
\newcommand{\Qed}{{\hfill$\Box$}}

{\vspace{16pt}\begin{em}\noindent {\bf #1.}}%
{\vspace{10pt}\end{em}}

{\vspace{16pt}\begin{em}\noindent {\bf Corollary.
\arabic{section}.\stepcounter{com}\arabic{com}}}%
{\vspace{10pt}\end{em}}

{\vspace{16pt}\begin{rm}\noindent{\bf Definition.
\arabic{section}.\stepcounter{com}\arabic{com}}}%
{\vspace{10pt}\end{rm}}

{\vspace{16pt}\begin{em}\noindent {\bf Lemma.
\arabic{section}.\stepcounter{com}\arabic{com}}}%
{\vspace{10pt}\end{em}}

{\vspace{16pt}\begin{rm}\noindent {\bf Example.
\arabic{section}.\stepcounter{com}\arabic{com}}}%
{\vspace{10pt}\end{rm}}

{\begin{rm}\noindent {\bf Proof : }}%
{\hfill{$\vrule height5pt width7pt depth2pt$}\vspace{12pt} \end{rm}}

\begin{document}

\baselineskip=1.2\baselineskip
\parskip=0.35\baselineskip

\title{
{Weighted Composition Operators of $C_0(X)$'s}
\thanks{This research  is  partially  
supported  by 
National Science Council of Taiwan, R.O.C.,
under the grant NSC 82-0208-M110-071.}%
}

\author{%
Jyh-Shyang Jeang and Ngai-Ching Wong\\
Department of Applied Mathematics\\
National Sun Yat-sen University\\
Kao-hsiung, 80424, Taiwan, R.O.C.\\
}

\date{} 

\maketitle

\begin{abstract}
In this paper, we prove that into isometries and disjointness  
preserving linear maps 
from $C_0(X)$ into $C_0(Y)$ are essentially  weighted composition  
operators $Tf = h\cdot f\circ\varphi$
for some continuous map $\varphi$ and some continuous
scalar-valued function $h$.
\end{abstract}

\section{Introduction.}

Let $X$ and $Y$ be locally compact Hausdorff spaces. 
Let $C_0(X)$ (resp. $C_0(Y)$) be
the Banach space of continuous scalar-valued (i.e.
real- or complex-valued) functions defined on 
$X$ (resp. $Y$) vanishing at infinity and
equipped with the supremum norm.
The classical Banach-Stone theorem  gives a description of surjective
isometries from $C_0(X)$ onto $C_0(Y)$.
They are all
{\em weighted composition operators} $Tf = h\cdot f\circ\varphi$ (i.e. 
$Tf(y) = h(y)f(\varphi(y))$, $\forall y \in Y$) for 
some homeomorphism $\varphi$ from
$Y$ onto $X$ and some continuous 
scalar-valued function $h$ on $Y$ with $|h(y)| \equiv 1$, $\forall  
y \in Y$. 
Different generalizations (see e.g. \cite{Beh88}, \cite{Cam70},  
\cite{Jar84}, \cite{Jar88},
\cite{Ves95}) of the Banach-Stone Theorem have been
studied in many years.  
Some of them discuss the structure of {\em into} isometries and  
disjointness preserving
linear maps (see e.g. \cite{Hol66}, \cite{Jar90}).
A linear map from $C_0(X)$ into $C_0(Y)$ is said to be 
{\em disjointness preserving} if $f\cdot g = 0$ in $C_0(X)$ implies
$Tf\cdot Tg = 0$ in $C_0(Y)$.
In this paper,
we shall discuss the structure of weighted composition operators from
$C_0(X)$ into $C_0(Y)$.
We prove that every into 
isometry and every disjointness preserving linear
map from $C_0(X)$ into $C_0(Y)$ is essentially a weighted composition
operator. 

{\bf Theorem 1.} {\em
Let $X$ and $Y$ be locally compact Hausdorff spaces and $T$ a linear
isometry from $C_0(X)$ into $C_0(Y)$.  Then there exist a locally\ 
compact subset $Y_1$ (i.e. $Y_1$ is locally compact in the subspace  
topology)
and a weighted composition operator $T_1$ from $C_0(X)$ into $C_0(Y_1)$
such that for all $f$ in $C_0(X)$,
$$ 
Tf_{|_{Y_1}} = T_1f = h\cdot f\circ\varphi,
$$
for some quotient map $\varphi$ from $Y_1$
onto $X$ and some continuous scalar-valued
function $h$ defined on $Y_1$ with $|h(y)|\ \equiv 1$, 
$\forall y \in Y_1$. 
}

{\bf Theorem 2.} {\em
Let $X$ and $Y$ be locally compact Hausdorff spaces and $T$  a 
bounded disjointness 
preserving linear map from $C_0(X)$ into $C_0(Y)$. 
Then there exist an open subset $Y_1$ of $Y$  
and a weighted composition operator $T_1$ from $C_0(X)$ 
into $C_0(Y_1)$ such that for all $f$ in $C_0(X)$, $Tf$ vanishes  
outside $Y_1$ and
$$
Tf_{|_{Y_1}} = T_1f = h\cdot f\circ\varphi, 
$$
for some continuous map $\varphi$ from $Y_1$
into $X$ and some continuous scalar-valued
function $h$ defined on $Y_1$ with $h(y) \neq 0$, 
$\forall y \in Y_1$. 
}

Since weighted composition operators from $C_0(X)$ into $C_0(Y)$  
are disjointness
preserving, Theorem 2 gives a complete description of all such maps.
When $X$ and $Y$ are both compact, Theorems 1 and 2 reduce to the  
results
of W. Holsztynski \cite{Hol66} and K. Jarosz \cite{Jar90}, respectively.
It is plausible to think that Theorems 1 and 2 could be easily obtained
from their compact space versions by simply extending an into isometry
(or a bounded disjointness preserving linear map)
$T : C_0(X) \lrw C_0(Y)$ to a bounded linear map $T_\infty :  
C(X_\infty) \lrw
C(Y_\infty)$ of the same type, where $X_\infty = X\bigcup\{\infty\}$ and
$Y_\infty = Y\bigcup\{\infty\}$ are the one-point compactifications of 
the locally compact Hausdorff spaces $X$
and $Y$, respectively.  However, the example given in Section 4  
will show  that this
idea is sometimes fruitless because $T$ can have {\em no} such 
extensions at all.  We thus have to modify, and in some cases give
new arguments to, the proofs of W. Holsztynski \cite{Hol66} and K.  
Jarosz
\cite{Jar90} to fit into our more general settings in this paper.

Recall that for $f$ in $C_0(X)$, the {\em cozero} 
of $f$ is  $\coz(f) = \{x\in X : f(x) \not= 0\}$ and 
the {\em support} $\supp(f)$ of  $f$ is the closure of $\coz(f)$ in 
$X_{\infty}$.  
A linear map $T : C_0(X) \lrw C_0(Y)$ is disjointness preserving if 
$T$ maps functions with disjoint cozeroes to functions with  
disjoint cozeroes.
For $x$ in $X$, $\delta_x$ denotes the point evaluation 
at $x$, that is, $\delta_x$ is the linear functional on $C_0(X)$
defined by $\delta_x(f) = f(x)$.
For $y$ in $Y$, let $\supp(\delta_y \circ T)$ be the set of all $x$ in 
$X_{\infty}$ such that for any open neighborhood $U$ of $x$ in  
$X_{\infty}$
there is an $f$ in $C_0(X)$ with $Tf(y) \not= 0$ and $\coz(f)  
\subset U$.
The kernel of a function $f$ is denoted by $\ker f$.

\section{Isometries from $C_0(X)$ into $C_0(Y)$.}

{\bf Definition.} 
Let $X$ and $Y$ be locally compact Hausdorff spaces.  A map $\varphi$
from $Y$ into $X$ is said to be {\em proper} if preimages of  
compact subsets of $X$
under $\varphi$ are compact in $Y$.

It is obvious that $\varphi$ is proper if and only if  
$\lim_{y\to\infty} \varphi(y) = \infty$.
As a consequence,  a  proper continuous map $\varphi$ from a  
locally compact Hausdorff space
$Y$ onto a locally compact Hausdorff space $X$ is a quotient map, i.e.
$\varphi^{-1}(O)$ is open in $Y$ if and only if $O$ is open in $X$.  
A quotient map from a locally
compact space onto another is, however, not
necessarily proper.  For example, the quotient map $\varphi$ from  
$(-\infty,+\infty)$ onto 
$[0,+\infty)$ defined by
$$
\varphi(y) =
\left\{
\begin{array}{ll}
y, & y > 0, \\
0, & y \leq 0
\end{array}\right.
$$
is not proper.

{\bf Lemma 3.}  {\em 
Let $X$ and $Y$ be locally compact Hausdorff spaces, $\varphi$ a  
map from $Y$ 
into $X$, and $h$ a continuous scalar-valued function defined on  
$Y$ with bounds $M, m > 0$ such
that $m \leq |h(y)| \leq M$, $\forall y \in Y$.  Then the weighted  
composition
$Tf = h\cdot f\circ\varphi$ defines a (necessarily bounded) linear  
map from $C_0(X)$ into $C_0(Y)$
if and only if $\varphi$ is continuous and proper.}

{\sc Proof.} 
For the sufficiency, we need to verify that $h\cdot f\circ\varphi$  
vanishes at $\infty$ for all
$f$ in $C_0(X)$.  For any $\epsilon > 0$, $|f(x)| < \epsilon/M$  
outside some compact subset $K$ of $X$.
Since $\varphi$ is proper, $\varphi^{-1}(K)$ is compact in $Y$.   
Now the fact that $|h(y)\cdot f(\varphi(y))|
\leq M|f(\varphi(y))| < \epsilon$ outside $\varphi^{-1}(K)$  
indicates that $h\cdot f\circ\varphi \in C_0(Y)$.
The boundedness of $T$ is trivial in this case.

For the necessity, we first check the continuity of $\varphi$.   
Suppose $y_\lambda \lrw y$
in $Y$.  We want to show that $x_\lambda = \varphi(y_\lambda) \lrw  
\varphi(y)$ in $X$.  
Suppose not, by passing to a subnet if necessary, we can assume  
that $x_\lambda$ either
converges to some $x \neq \varphi(y)$ in $X$ or $\infty$.  If  
$x_\lambda \lrw x$ in $X$
then for all $f$ in $C_0(X)$,
\begin{eqnarray*}
h(y)f(x) &=& \lim h(y_\lambda)f(x_\lambda) = \lim  
h(y_\lambda)f(\varphi(y_\lambda))\\
         &=&  \lim Tf(y_\lambda) = Tf(y) = h(y)f(\varphi(y)).
\end{eqnarray*}
As $h(y) \neq 0$, $f(x) = f(\varphi(y))$, $\forall f \in C_0(X)$.   
Consequently,
we obtain a contradiction $x = \varphi(y)$.  If $x_\lambda \lrw  
\infty$ then
a similar argument gives $f(\varphi(y)) = 0$ for all $f$ in  
$C_0(X)$.  Hence
$\varphi(y) = \infty$, a contradiction again.  Therefore, $\varphi$  
is continuous from
$Y$ into $X$.  Finally, let $K$ be a compact subset of $X$ and we  
are going to see
that $\varphi^{-1}(K)$ is compact in $Y$, or equivalently, closed  
in $Y_\infty = Y \cup\{\infty\}$, the
one-point compactification of $Y$.  To see this, suppose $y_\lambda  
\lrw y$ in $Y_\infty$ and
$x_\lambda = \varphi(y_\lambda) \in K$.  We want $y \in  
\varphi^{-1}(K)$, i.e. $y \neq \infty$
and $\varphi(y) \in K$.  Without loss of generality, we can assume  
that $x_\lambda \lrw x$ for
some $x$ in $K$.  Now,
$$
\lim |Tf(y_\lambda)| = \lim |h(y_\lambda)f(\varphi(y_\lambda))|  
\geq m\lim |f(x_\lambda)| = m|f(x)|
$$
for all $f$ in $C_0(X)$.  This implies that $y \neq \infty$ and  
then a similar argument gives
$\varphi(y) = x \in K$. \Qed

The assumption on the bounds of $f$ in Lemma 3 is significant.
For example, let $X = Y = {\Bbb R} = (-\infty,+\infty)$ and define
$$
h(y)=\left\{
\begin{array}{ll}
e^y, & y<0,\\
1, & y \geq 0,
\end{array}\right.
\mbox{\qquad and \qquad}
\varphi(y)=\left\{
\begin{array}{ll}
\sin y, & y< 0,\\
y,& y \geq 0.
\end{array}\right.
$$
Then the weighted composition operator $Tf = h\cdot f\circ\varphi$ 
from $C_0({\Bbb R})$ into $C_0({\Bbb R})$ is well-defined.  
It is not difficult to see that ${\varphi}^{-1}([-{1\over  
2},{1\over 2}])$ is
not compact in ${\Bbb R}$.
On the other hand, if we redefine $h(y) = e^y$ and $\varphi(y) = y$  
for all $y$ in ${\Bbb R}$
then the weighted composition operator $T$ is not well-defined 
from $C_0({\Bbb R})$ into $C_0({\Bbb R})$, even though $\varphi$ is  
proper and continuous in this case.

Recall that a bounded linear map $T$ from a Banach space $E$ into a  
Banach space $F$ is
called an {\em injection} if there is an $m > 0$ such that $\|Tx\|  
\geq m\|x\|$, $\forall x \in E$.
It follows from the open mapping theorem that $T$ is an injection  
if and only if $T$ is one-to-one
and has closed range.

{\bf Proposition 4.} {\em
Let $X$ and $Y$ be locally compact Hausdorff spaces,  
$\varphi$ a map from $Y$ into
$X$ and $h$ a continuous scalar-valued function defined on $Y$.
The weighted composition operator $Tf = h\cdot f\circ\varphi$
from $C_0(X)$ into $C_0(Y)$ is an injection if and only if  
$\varphi$ is continuous, proper and onto
and $h$ has bounds $M, m > 0$ such that
$m \leq |h(y)| \leq M$, $\forall y \in Y$.  In this case, $\varphi$  
is a quotient map and thus $X$ is a quotient space
of $Y$.}

{\sc Proof.}
The sufficiency follows easily from Lemma 3 and the observation  
that $\|Tf\| = \|h\cdot f\circ\varphi\|
\geq m\|f\|$, $\forall f \in C_0(X)$.
For the necessity, we first note that there are constants $M, m >  
0$ such that $m\|f\| \leq \|Tf\| \leq M\|f\|$
for all $f$ in $C_0(X)$.  It is then obvious that $m \leq |h(y)|  
\leq M$, $\forall y \in Y$.  By Lemma 3, $\varphi$
is continuous and proper.  Finally, we check that $\varphi$ is  
onto.  It is not difficult to see that $\varphi$
has dense range.  In fact, if $\varphi(Y)$ were not dense in $X$,  
then there were an $x$ in $X$ and a
neighborhood $U$ of $x$ in $X$ such that 
$U\cap\varphi(Y) = \emptyset$.  Choose an $f$ in $C_0(X)$ such that
$f(x) = 1$ and $f$ vanishes outside $U$.  Then $Tf(y) =
h(y)f(\varphi(y)) = 0$ for all $y$ in $Y$, i.e. $Tf = 0$.  Since  
$T$ is an injection, we get a contradiction that
$f = 0$.  We now show that $\varphi(Y) = X$.  Let $x \in X$ and $K$  
a compact neighborhood of $x$ in $X$. 
By the density of $\varphi(Y)$ in $X$, there is a net  
$\{y_\lambda\}$ in $Y$ such that $x_\lambda
= \varphi(y_\lambda) \lrw x$ in $X$.  Without loss of generality,  
we can assume that $x_\lambda$ belongs to
$K$ for all $\lambda$.  Since $\varphi^{-1}(K)$ is compact in $Y$,  
$\varphi(\varphi^{-1}(K))$ is a compact
subset of $X$ containing the net $\{x_\lambda\}$.  Consequently, $x  
= \lim x_\lambda$ belongs to 
$\varphi(\varphi^{-1}(K)) \subset \varphi(Y)$.  \Qed

{\sc Proof of Theorem 1.}
We adopt some notations from W. Holsztynski \cite{Hol66} and K.  
Jarosz \cite{Jar90}.  
Let $X_\infty = X\cup\{\infty\}$ and $Y_\infty = Y\cup\{\infty\}$  
be the one-point compactifications of
$X$ and $Y$, respectively.
For each $x$ in $X$ and $y$ in $Y$, put
\begin{eqnarray*}
S_{x} &=& \{ f \in C_0(X) :  |f(x)| = \|f\| = 1\},\\
R_{y} &=& \{ g \in C_{0}(Y) : |g(y)| = \|g\| = 1\}, \mbox{ and }\\ 
Q_{x} &=& \{ y \in Y : T(S_{x}) \subset R_{y}\}.
\end{eqnarray*}

We first claim that $\{Q_{x}\}_{x\in X}$ is a disjoint family of  
non-empty subsets of $Y$.  
In fact,
for $f_1,\ f_2,\cdots,\ f_n$ in $S_x$,  
let $h = \sum_{i=1}^{n}\overline{f_i(x)} f_i$. Then $\|h\|=n$
and thus $\|Th\| = n$.  Hence there is a $y$ in $Y$ such that 
$|\sum_{i=1}^{n} \overline{f_i(x)} Tf_i(y)| = |Th(y)| = n$.  This  
implies $|Tf_i(y)| = 1$ for all $i=1, 2,\cdots, n$.
In other words, $y \in \bigcap_{i=1}^{n} (Tf_i)^{-1}(\Gamma)$, where
$\Gamma =\{ z : |z| = 1\}$.   We have just proved that
the family $\{ (Tf)^{-1}(\Gamma) : f \in S_x \}$ of closed subsets  
of the
compact space $Y_{\infty}$ has finite intersection
property.  
It is plain that $\infty \notin (Tf)^{-1}(\Gamma)$ for all $f$ in $S_x$.
Hence  $Q_{x} = \bigcap_{f \in S_x}(Tf)^{-1}(\Gamma)$ is non-empty  
for all $x$ in $X$.  
Moreover, $Q_{x_1}\cap Q_{x_2} = \emptyset$ if $x_1 \neq x_2$ in  
$X$.  In fact, $f_1$ in $S_{x_1}$ and $f_2$ in $S_{x_2}$ exist such  
that
$\coz(f_1)\cap\coz(f_2) = \emptyset$.  If 
there is a $y$ in $Q_{x_1}\cap Q_{x_2}$ then it follows from $Tf_1  
\in R_y$ and $Tf_2 \in R_y$ that
$1 = \|f_1 + f_2\| = \|T(f_1 + f_2)\| = |T(f_1 + f_2)(y)| = 2$, a  
contradiction.

Let
$Y_1 = \bigcup_{x\in X} Q_x$.  It is
not difficult to see that $\supp(\delta_{y} \circ T)$ = $\{x\}$
whenever $y\in Q_x$.  So we can define a surjective map $\varphi : Y_{1}
\rightarrow X$ by 
$$
\{\varphi (y)\} = \supp(\delta_y \circ T).
$$  
Note that for all $f$ in $C_0(X)$ and for all $y$ 
in $Y_{1}$,
\begin{equation}
\varphi(y) \notin \supp(f) \Longrightarrow T(f)(y) =0. 
\end{equation}
In fact, if $Tf(y) \neq 0$, without loss of generality, we can assume
$Tf(y) = r > 0$ and $\|f\| = 1$.  Since $\varphi(y) \notin \supp(f)$,
there is a $g$ in $C_0(X)$ such that $\coz(f) \bigcap \coz(g) =
\emptyset$ and $Tg(y) = \|g\| = 1$.  Hence $1+r = T(f+g)(y) >  
\|f+g\| =1$,
a contradiction.

Now, we want to show that $\varphi$ is continuous.  Suppose  
$\varphi$ were
not continuous at some $y$ in $Y_{1}$, without loss of generality,
let $\{y_{\lambda}\}$ be a net converging to $y$ in $Y_1$ such that
$\varphi (y_{\lambda}) \rightarrow x \neq
\varphi (y)$ in $X_{\infty}$. 
Then there exist disjoint
neighborhoods $U_1$ and $U_2$ of $x$ and $\varphi (y)$ in  
$X_\infty$, respectively, and 
a $\lambda_{0}$ such that $\varphi (y_{\lambda}) \in U_{1}$, $\forall
\lambda \geq \lambda_{0}$.  Let $f \in C_0(X)$ such that
coz($f) \subseteq U_{2}$ and $T(f)(y) = \| f \| = 1$.  As $\supp(f)  
\bigcap
U_{1} = \emptyset$, we have $\varphi (y_{\lambda}) \notin \supp(f)$,
  $\forall \lambda \geq \lambda_{0}$. By (1), $T(f)(y_{\lambda}) = 0$,
$\forall \lambda \geq \lambda_{0}$.  This implies $T(f)$ is not  
continuous
at $y$, a contradiction.

For each $y$ in $Y_{1}$, put
\begin{eqnarray*}
J_{y} &=& \{ f \in C_0(X) : \varphi (y) \notin \supp(f) \}, \mbox{  
and}\\
K_{y} &=&  \{ f \in C_0(X) : f(\varphi (y))=0 \}.
\end{eqnarray*}
For $f$ in $K_y$ and $\varepsilon > 0$, let $X_1 = \{x\in X :  
|f(x)| \geq
\varepsilon\}$ and $X_2 = \{x\in X : |f(x)| \leq \varepsilon/2\}$.  Let
$g$ be a continuous function defined on $X$ such that $0 \leq g(x)  
\leq 1, \forall x \in X$,
$g(x) = 1, \forall x\in X_1$, and $g(x) = 0, \forall x \in X_2$.
Let $f_{\varepsilon} = g\cdot f$.  Then $f_{\varepsilon} \in J_y$ and 
$\|f_{\varepsilon}-f\| \leq 2\varepsilon$.
One thus can show that $J_{y}$ is a dense subset of $K_{y}$.  By (1),
$J_{y} \subset \ker(\delta_{y} \circ T)$, and hence
${\ker}(\delta_{ \varphi(y)}) = K_y \subset {\ker}(\delta_{y} \circ T)$.
Consequently, there exists a scalar $h(y)$ such that 
$\delta_{y} \circ T = h(y) \cdot \delta_{ \varphi(y)}$, i.e.
$$
T(f)(y) = h(y) \cdot f(\varphi(y)), \qquad \forall f \in C_0(X).
$$
It follows from the definition of $Y_1$ that $h$ is continuous
on $Y_1$ and  $|h(y)|=1$, $\forall y \in Y_1$.

It is the time to see that $Y_{1}$ is locally compact. 
For each $y_1$ in $Y_1$ and a neighborhood $U_1$ of $y_1$ in $Y_1$,  
we want to find
a compact neighborhood $K_1$ of $y_1$ in $Y_1$ such that $y_1\in  
K_1 \subset U_1$.
Let $x_1 = \varphi(y_1)$ in $X$.  Then 
$$
|Tf(y_1)| = |f(x_1)|, \quad \forall f \in C_0(X).
$$  
Fix $f_1$ in $S_{x_1}$.  Then $V_1 = \varphi^{-1}(\{ x \in X :  
|f_1(x)| > {1\over 2}\})\cap U_1$ 
is an open neighborhood of $y_1$ in $Y_1$ and contained in $U_1$.
Since $V_1 = W \bigcap Y_1$ for some neighborhood $W$ of $y_1$ in  
$Y$, there
exists a compact neighborhood $K$ of $y_1$ in $Y$ such that $y_1  
\in K \subset
W$.  We are going to verify that
$K_1 = K \bigcap Y_1$ is a compact neighborhood of $y_1$ in $Y_1$.
Let $\{y_{\lambda}\}$ be a net in $K_1 \subset V_1$.  By passing to  
a subnet, we can assume
that $y_{\lambda}$ converges to $y$ in $K$ and we want to show $y  
\in Y_1$.
Let $x_{\lambda} = \varphi(y_{\lambda})$ in $X$.  Since  
$X_{\infty}$ is compact,
by passing to a subnet again, we can assume that $x_{\lambda}$  
converges to $x$ in $X$ or 
$x_{\lambda} \rightarrow \infty$.  If $x_{\lambda} \rightarrow x$  
in $X$, $|Tf(y)| =
\lim|Tf(y_{\lambda})| = 
\lim |h(y_\lambda)f(\varphi(y_{\lambda}))| = \lim |f(x_{\lambda})|  
= |f(x)|$, for all $f$ in $C_0(X)$.
Hence $y \in Q_x$, and thus $y \in Y_1$.  If $x_{\lambda}  
\rightarrow \infty$,
$|Tf_1(y)| = \lim |Tf_1(y_{\lambda})| =
\lim |h(y_\lambda)f_1(\varphi(y_{\lambda}))| = \lim  
|f_1(x_{\lambda})| = 0$.
However, the fact that $y_{\lambda} \in V_1$ ensures  
$|Tf_1(y_{\lambda})| = 
|f_1(x_{\lambda})| > {1 \over 2}$ for all $\lambda$, a contradiction.
Hence $Y_1$ is locally compact.

Let $T_1 : C_0(X) \rightarrow C_0(Y_1)$ defined by $T_1 f = h \cdot  
f\circ\varphi
$.  It is clear that $T_1$ is a linear isometry and $Tf_{|_{Y_1}} =  
T_1 f$.  By Proposition 4,
the surjective continuous map $\varphi$ is proper and thus a  
quotient map.
The proof is complete.
\Qed

In Theorem 1, $Y_{1}$ can be neither open nor closed in $Y$ and  
$\varphi$ may not be
an open map.  See the following examples.

{\bf Example 5.}
Let $X=[0,+\infty)$ and $Y=[-\infty,+\infty]$.
Let $T$ be a linear isometry from $C_{0}(X)$ into $C_{0}(Y)$ defined
for all $f$ in $C_{0}(X)$ by
$$
Tf(y)=\left\{
\begin{array}{ll}
f(y), & 0 \leq y < +\infty,\\
{e^{y} \over 2}(f(-y)+f(0)), & -\infty < y \leq 0,\\
0, & y = \pm\infty. 
\end{array}\right.
$$  
Then in notations of Theorem 1,
$Y_{1}=[0,+\infty)$ is neither closed  nor open in $Y$.  In
this case,
$\varphi (y)=y$ for all $y$ in $[0,+\infty)$, and
$X$ and $Y_1$ are homeomorphic. \Qed

{\bf Example 6.}   
Let $X = {\Bbb R}$ and 
$Y = \{(x,y) \in {\Bbb R}^2 : y = 0\} 
\bigcup \{(x,y) \in {\Bbb R}^2 : 0 \leq x, 0 \leq y \leq 1\}$.
Let $\varphi : Y \rightarrow X$ defined by $\varphi(u_1,u_2) =  
u_1$.  Then
$\varphi$ is continuous, onto and proper, and thus a quotient map.  
Moreover, $T : C_0(X) \rightarrow C_0(Y)$ defined
by $Tf = f\circ \varphi$ is a linear isometry.  
Note that $O = \{(x,y) \in {\Bbb R}^2 : 0 \leq x < 1, 0 < y \leq 1\}$ 
is open in $Y$,
but $\varphi(O) = [0,1)$ is not open in $X$. Hence $\varphi$ is not an 
open map. \Qed

\section{Disjointness preserving linear maps from $C_0(X)$ to $C_0(Y)$.}

It is clear that Theorem 2 follows from the following more general  
result in which
discontinuity of the linear disjointness preserving map $T$ is  
allowed.  The payoff
of the discontinuity is a finite subset $F$ of $X$ at which the  
behaviour of $T$ is
not under control.

{\bf Theorem 7.} {\em
Let $X$ and $Y$ be locally compact Hausdorff spaces and
$T$ a disjointness preserving linear map from $C_0(X)$ into $C_0(Y)$.
Then $Y$ can be written as a disjoint union $Y = Y_1 \cup Y_2 \cup  
Y_3$, in which
$Y_2$ is open and $Y_3$ is closed.  A continuous
map $\varphi$ from $Y_1 \cup Y_2$ into $X_{\infty}$ 
exists such that for every $f$ in $C_0(X)$,
\begin{equation}
\varphi(y) \notin \supp(f) \Longrightarrow T(f)(y) =0. 
\end{equation}
Moreover,
a continuous bounded non-vanishing scalar-valued function $h$ on  
$Y_1$ exists 
such that 
\begin{eqnarray*}
Tf_{|_{Y_1}} &=& h\cdot f\circ\varphi, \mbox{ and}\\
Tf_{|_{Y_3}} &=& 0. \\
\end{eqnarray*}
Furthermore, $F = \varphi(Y_2)$ is a finite set
and the functionals $\delta_y \circ T$
are discontinuous on $C_0(X)$ for all $y$ in $Y_2$.}

{\sc Proof.}
We shall follow the plan of K. Jarosz in his compact space version  
\cite{Jar90}.
Set
\begin{eqnarray*}
Y_3 &=& \{ y\in Y |\delta_y \circ T \equiv 0 \},\\
Y_2 &=& \{ y\in Y | \delta_y \circ T \mbox{ is discontinuous} \},  
\mbox{ and}\\
Y_1 &=& Y \setminus ( Y_2 \bigcup Y_3 ).
\end{eqnarray*}

First, we claim that 
$\supp(\delta_y\circ T)$ contains exactly one point for every $y$  
in $Y_1\bigcup Y_2$.
Suppose on the contrary that $\supp(\delta_y\circ T)$ contains two  
distinct points $x_1$ and
$x_2$ in $X_\infty$.  Let $U_1$ and $U_2$ be neighborhoods of $x_1$  
and $x_2$ in $X_\infty$, 
respectively, such that $U_1\cap U_2 = \emptyset$.  Let $f_1$ and  
$f_2$ in $C_0(X)$ with
$\coz(f_1) \subset U_1$ and $\coz(f_2) \subset U_2$ be such that  
$Tf_1(y) \neq 0$ and $Tf_2(y) \neq 0$.
However, $f_1f_2 = 0$ implies $Tf_1 Tf_2 = 0$, a contradiction.
Suppose $\supp(\delta_y\circ T)$ is empty.  Then we can write the  
compact
Hausdorff space $X_\infty$ as a finite union of open sets $X_\infty  
= \bigcup_{i=1}^n U_i$
such that $Tf(y) = 0$ whenever $\coz(f) \subset U_i$ for some $i =  
1, 2, \ldots, n$.
Let ${\mbox{\boldmath $1$}} = \sum_{i=1}^n f_i$ be a continuous  
decomposition of the
identity subordinate to $\{U_i\}_{i=1}^n$.  Then for all $f$ in  
$C_0(X)$,
$Tf(y) = \sum_{i=1}^n T(ff_i)(y) = 0$.  This says $\delta_y\circ T  
\equiv 0$ and
thus $y \in Y_3$.

Next we define a map $\varphi$ from $Y_1\bigcup Y_2$ into 
$X_{\infty}$ by  
$$
\{ \varphi(y) \} = \supp(\delta_y\circ T).
$$
We now prove (2).  Assume $\varphi(y) \notin \supp(f)$.  Then there  
is an open neighborhood
$U$ of $\varphi(y)$ disjoint from $\coz(f)$.  Let $g \in C_0(X)$  
such that $\coz(g) \subset U$ and 
$Tg(y) \neq 0$.  Since $fg = 0$ and $T$ is disjointness preserving,  
$Tf(y) = 0$ as asserted.  

It then follows from (2) the continuity of $\varphi$ as one can  
easily modify an argument of the proof
of Theorem 1 for this goal.  Similarly, it also follows from (2)  
the desired representation
\begin{equation}
Tf(y) = h(y)f(\varphi(y)), \qquad \forall f \in C_0(X), \forall y  
\in Y_1,
\end{equation}
where $h$ is a continuous non-vanishing scalar-valued function  
defined on $Y_1$.

{\bf Claim.}  {\em Let $\{y_n\}_{n=1}^{\infty}$ be a sequence in  
$Y_1\bigcup 
Y_2$ such that $x_n = \varphi(y_n)$'s are distinct points of $X$.  Then
$$
\limsup\|\delta_{y_n}\circ T\| < \infty.
$$
In particular, only finitely many $\delta_y\circ T$ can have  
infinite norms.}

Assume the contrary and, by passing to a subsequence if necessary,  
we have
$$\|\delta_{y_n}\circ T\| > n^4,\ \ \ n = 1,2,\cdots.$$
Let $f_n \in C_0(X)$ with $\|f_n\| \leq 1$ such that  
$$|Tf_n(y_n)|\geq n^3,\ 
n = 1,2,\cdots.$$  Let $V_n$, $W_n$ and $U_n$ be open subsets of  
$X$ such that
$x_n \in V_n \subseteq \overline{V_n} \subseteq W_n \subseteq  
\overline{W_n}
\subseteq U_n$ and $U_n \bigcap U_m = \emptyset$ if $n \not= m,\  
n,m = 1,2,
\cdots$, and let $g_n \in C(X_{\infty})$ such that $0 \leq g_n \leq 1,\ 
{g_n}_{|_{V_n}} \equiv 1$ and ${g_n}_{|_{X_{\infty}\setminus W_n}}  
\equiv 0$,
$n = 1,2,\cdots$.  Then (2) implies
$$\begin{array}{ccl}
Tf_n(y_n) & = & T(f_n g_n)(y_n) + T(f_n(1-g_n))(y_n)\\
 & = & T(f_n g_n)(y_n),\ \ \ n = 1,2,\cdots.
\end{array}$$
Therefore, we can assume $\supp f_n \subset U_n$.  Let $f =  
\sum_{n=1}^{\infty}
{1\over{n^2}} f_n$  in $C_0(X)$.  By (2) again,
$|Tf(y_n)| = | {1\over{n^2}} Tf_n(y_n) | \geq n$ for $n = 1,2,\cdots$.
This conflicts with the boundedness of $Tf$ in $C_0(Y)$, and the  
claim is thus
verified.

The assertion $F = \varphi(Y_2)$ is a finite subset of $X$ is  
clearly a consequence of
the claim while the boundedness of  $h$ follows from the claim and (3). 
It is also plain that $Y_3 = \bigcap\{\ker Tf : f \in C_0(X)\}$ is  
closed in $Y$.  Finally, to
see that $Y_2$ is open, we consider for every $f$ in $C_0(X)$,
\begin{eqnarray*}
	\sup \{|Tf(y)| : y \in \overline{Y_1\cup Y_3}\} 
&=& 	\sup \{|Tf(y)| : y \in {Y_1\cup Y_3}\}\\
&=& 	\sup \{|Tf(y)| : y \in {Y_1}\}\\
&=& 	\sup \{|h(y)f(\varphi(y))| : y \in {Y_1}\} \leq M\|f\|,
\end{eqnarray*}
where $M > 0$ is a bound of $h$ on $Y_1$.  It follows that the  
linear functional 
$\delta_y\circ T$ is bounded for all $y$ in $\overline{Y_1\cup  
Y_3}$, and thus
$Y_2\cap\overline{Y_1\cup Y_3} = \emptyset$.  Hence, $Y_1\cup Y_3 =  
\overline{Y_1\cup Y_3}$
is closed.  In other words, $Y_2$ is open.
\Qed

{\bf Theorem 8.} {\em
Let $X$ and $Y$ be locally compact Hausdorff spaces and $T$ 
a bijective disjointness preserving linear map from $C_0(X)$ onto
$C_0(Y)$.  Then $T$ is a bounded weighted composition operator, and
$X$ and $Y$ are homeomorphic.}

{\sc Proof.}
We adopt the notations used in Theorem 7.  Since $T$ is surjective,  
$Y_3 = 
\emptyset$.  We are going to verify that $Y_2 = \emptyset$, too.  First,
we note that the finite set $F\setminus\{\infty\}$ consists of  
non-isolated 
points in $X$.
In fact, if $y\in Y_2$ such that $x = \varphi(y)$ is an isolated  
point in $X$
then it follows from (2) that for every $f$ in $C_0(X)$, $f(x) = 0$  
implies
$\varphi(y) = x \not\in \supp f$ and thus $Tf(y) = 0$.  Hence,  
$\delta_y \circ
T = \lambda \delta_x$ for some scalar $\lambda$.  Therefore,  
$\delta_y \circ T$
is continuous, a contradiction to the assumption that $y\in Y_2$.   
We then
claim that $\varphi(Y) = \varphi(Y_1\bigcup Y_2)$ is dense in $X$.   
In fact,
if a nonzero $f$ in $C_0(X)$ exists such that $\supp f\bigcap  
\varphi(Y) =
\emptyset$ then $ Tf = 0$ by (2), conflicting with the injectivity  
of $T$.
Since 
$$
X = \overline{\varphi(Y)}  =   
\overline{\varphi(Y_1)\bigcup\varphi(Y_2)}  
=  \overline{\varphi(Y_1)\bigcup F} =  \overline{\varphi(Y_1)}  
\mbox{ or } \overline{\varphi(Y_1)}\bigcup\{\infty\},
$$
for every $f$ in $C_0(X)$,
$$
Tf_{|_{Y_1}} = 0  \Longrightarrow  f_{|_{\varphi(Y_1)}} = 0 
 \Longrightarrow  f = 0
 \Longrightarrow  Tf_{|_{Y_2}} = 0.
$$
Therefore, the open set $Y_2 = \emptyset$ by the surjectivity of $T$.
Theorem 7 then gives 
$$Tf = h\cdot (f\circ\varphi),\ \ \ \forall f\in C_0(X).$$
This representation implies that $T^{-1}$ is also a bijective  
disjointness 
preserving linear map from $C_0(Y)$ onto $C_0(X)$.  The above  
discussion 
provides that 
$$T^{-1}g = h_1 \cdot g\circ\varphi_1,\ \ \ \forall g\in C_0(Y),$$
for some continuous non-vanishing scalar-valued function $h_1$ on $X$
and continuous function $\varphi_1$ from $X$ into $Y$.  It is plain that
$\varphi_1 = \varphi^{-1}$ and thus $X$ and $Y$ are homeomorphic.
\Qed

\section{A counter example.}

The following example shows that not every into isometry or bounded  
disjointness
preserving linear map from
$C_0(X)$ into $C_0(Y)$ can be extended to a bounded linear map 
from $C(X_\infty)$ into $C(Y_\infty)$ of the same type.  Here $X$  
and $Y$ are
locally compact Hausdorff spaces with one-point compactifications  
$X_\infty$ and
$Y_\infty$, respectively.

{\bf Example 9.}
Let $X = [0,+\infty)$, $Y = (-\infty,+\infty)$ and 
the underlying scalar field is the field ${\Bbb R}$ of real numbers. 
Let
$$
h(y) = \left\{
\begin{array}{cl}
1, & y > 2, \\
y-1, & 0\leq y \leq 2,\\ 
-1, & y < 0,
\end{array}\right.
$$
and
$$
\varphi(y)= \left\{
\begin{array}{rl}
y,& y \geq 0,\\
-y,& y < 0.   
\end{array}\right.
$$
Then the weighted composition operator
$Tf = h \cdot f \circ \varphi$ is simultaneously an into isometry and 
a bounded disjointness preserving
linear map from $C_0([0,+\infty))$ into $C_0((-\infty,+\infty))$.   
However,
no bounded linear extension $T_\infty$ from 
$C([0,\infty])$ into $C((-\infty,+\infty)\cup\{\infty\})$ of
$T$ can be an into isometry or a disjointness preserving linear map. 

Suppose, on the contrary, $T_\infty$ were an into isometry.  
Consider $f_n$ in 
$C_0([0,+\infty))$ defined by
$$
f_n(x)= \left\{
\begin{array}{cl}
1, & 0 \leq x \leq n,\\
\frac{2n-x}{n}, & n < x < 2n, \quad n = 1, 2, \ldots.\\
0, & 2n \leq x \leq +\infty,    
\end{array}\right.
$$
Note that $\delta_y\circ T_\infty$ can be considered as a bounded Borel
measure $m_y$ on $[0,+\infty]$ for all point evaluation $\delta_y$ at
$y$ in $(-\infty,+\infty)\cup\{\infty\}$ with total variation 
$\|m_y\| = \|\delta_y\circ T_\infty\| \leq 1$.  Let  
$\mbox{\boldmath $1$}$ be
the constant function $\mbox{\boldmath $1$}(x) \equiv 1$ in  
$C([0,+\infty])$. 
For all $y$ in $(-\infty,+\infty)$,
\begin{eqnarray*}
T_\infty\mbox{\boldmath $1$}(y) &=& \delta_y\circ  
T_\infty(\mbox{\boldmath $1$})
= \int_{[0,+\infty]} \mbox{\boldmath $1$}\, dm_y \\
 &=& \lim_{n\to \infty} \int_{[0,+\infty]} f_n\,dm_y + m_y(\{\infty\})
= \lim_{n\to \infty} \delta_y\circ T_\infty(f_n) + m_y(\{\infty\}) \\
 &=& \lim_{n\to \infty} Tf_n(y) + m_y(\{\infty\})
= \lim_{n\to \infty} h(y)\cdot f_n(\varphi(y)) + m_y(\{\infty\}) \\
 &=&  h(y) + m_y(\{\infty\}).  
\end{eqnarray*} 
Let $g(y) = m_y(\{\infty\})$ for all $y$ in $(-\infty,+\infty)$.  Then
$g(y) = T_\infty\mbox{\boldmath $1$}(y) - h(y)$ is continuous on  
$(-\infty,+\infty)$ 
and $|g(y)| = |m_y(\{\infty\})| \leq \|m_y\| \leq 1$, $\forall y  
\in (-\infty,+\infty)$.
Note that $\|T_\infty\mbox{\boldmath $1$}\| = 1$.  Therefore, $g(y)  
= T_\infty\mbox{\boldmath $1$}(y) - 1
\leq 0$ when $y > 2$, and $g(y) = T_\infty\mbox{\boldmath $1$}(y) +  
1 \geq 0$ when
$y < -2$.  We claim that $g(y)g(-y) = 0$ whenever $|y| > 2$.  
In fact, if for example 
$g(y_0) < -\delta$ for some $y_0 > 2$ and some $\delta > 0$,
then for each small $\epsilon > 0$, $0 \leq T_\infty\mbox{\boldmath  
$1$}(y) < 1 - \delta$
for all $y$ in $(y_0-\epsilon,y_0+\epsilon)$.  We can choose an $f$  
in $C_0([0,+\infty))$
satisfying that $f(y_0) = \|f\| = 1$ and $f$ vanishes outside  
$(y_0-\epsilon,y_0+\epsilon)
\subset (2,+\infty)$.  Now,
\begin{eqnarray*}
T_\infty(\mbox{\boldmath $1$} + \delta f)(y) 
&=& T_\infty(\mbox{\boldmath $1$})(y) + \delta T_\infty(f)(y)\\
&=& T_\infty(\mbox{\boldmath $1$})(y) + \delta T(f)(y)\\
&=& h(y) + g(y) +\delta h(y)f(\varphi(y))\\
&=& \left\{\begin{array}{cc}
1 + g(y) + \delta f(y), & y > 2, \\
T_\infty\mbox{\boldmath $1$}(y), & -2\leq y \leq 2,\\ 
-1 + g(y) - \delta f(-y), & y < -2.
\end{array}\right.
\end{eqnarray*}
Since $\|T_\infty(\mbox{\boldmath $1$} + \delta f)\| = \|1 + \delta  
f\| = 1 + \delta$ and
$|T_\infty(\mbox{\boldmath $1$} + \delta f)(y)| \leq 1$ unless $-y  
\in (y_0-\epsilon,y_0+\epsilon)$,
there is a $y_1$ in $(y_0-\epsilon,y_0+\epsilon)$ such that
$|-1 + g(-y_1) - \delta f(y_1)| = 1 + \delta$.  It forces that  
$g(-y_1) = 0$.  Since $\epsilon$ can be
arbitrary small, we have $g(-y_0) = 0$ and our claim that  
$g(y)g(-y) = 0$ whenever $|y| > 2$ 
has thus been verified.
As $T_\infty\mbox{\boldmath $1$}$ is continuous on  
$(-\infty,+\infty)\cup\{\infty\}$, we must have
$$
\lim_{y\to+\infty} T_\infty\mbox{\boldmath $1$}(y) =  
\lim_{y\to-\infty} T_\infty\mbox{\boldmath $1$}(y),
$$
that is,
$$
\lim_{y\to+\infty} -1 + g(y) = \lim_{y\to-\infty} 1 + g(y).
$$
Let $L$ be their common (finite) limit.  Then 
$$
\lim_{y\to+\infty} g(y) = L + 1, \qquad \qquad \lim_{y\to-\infty}  
g(y) = L - 1.
$$
Consequently,
$$
0 = \lim_{y\to+\infty} g(y)g(-y) = L^2 - 1.
$$
It follows that $L = \pm 1$, and thus either $\lim_{y\to+\infty}  
g(y) = 2$ or 
$\lim_{y\to-\infty} g(y) = -2$.  Both of them contradicts the fact that
$|g(y)| \leq 1$, $\forall y \in (-\infty,+\infty)$.

	On the other hand, suppose $T_\infty$ were disjointness  
preserving.  Since
$f_n(\mbox{\boldmath $1$} - f_{2n}) = 0$, we have $T_\infty  
f_n\cdot T_\infty(\mbox{\boldmath $1$} - f_{2n}) = 0$.
That is,
$$
T_\infty f_n(y)\cdot T_\infty(\mbox{\boldmath $1$} - f_{2n})(y) = 0, 
\quad \forall y \in (-\infty,+\infty)\cup\{\infty\}.
$$
When $|y| < n$ and $y \neq 1$, $T_\infty f_n(y) = Tf_n(y) = h(y)  
\neq 0$ and
hence $T_\infty(\mbox{\boldmath $1$})(y) = T_\infty(f_{2n})(y) = 
T(f_{2n})(y) = h(y)$.   Since 
$T_\infty\mbox{\boldmath $1$}$ is continuous on  
$(-\infty,+\infty)\cup\{\infty\}$, we must have
$$
+1 = \lim_{y\to+\infty} h(y) = \lim_{y\to-\infty} h(y) = -1,
$$
a contradiction again. \Qed

\bibliographystyle{plain}

\end{document}